%
\magnification 1200
\input amssym.def
\nopagenumbers
\headline={\ifnum\pageno=1 \hfill \else\hss{\tenrm--\folio--}\hss \fi}
\font\titlefont=cmr10 scaled \magstep2
\font\mathtitlefont=cmmi10 scaled \magstep2
\font\mathsymtitlefont=cmsy10 scaled \magstep2
\font\addressfont=cmr10 at 10truept
\font\ttaddressfont=cmtt10 at 10truept
\let\sPP=\smallbreak
\let\mPP=\medbreak
\let\bPP=\bigbreak

\def\sLP{\smallbreak\noindent}
\def\mLP{\medbreak\noindent}
\def\bLP{\bigbreak\noindent}
\def\wrt{with respect to }
\def\st{such that}
\def\iy{\infty}
\def\const{{\rm const.}\,}
\def\LHS{left hand side}
\def\RHS{right hand side}
\def\half{{\scriptstyle{1\over2}}}
\def\thalf{{\textstyle{1\over2}}}
\def\rep{representation}
\let\pa=\partial
\let\lan=\langle
\let\ran=\rangle

\def\al{\alpha}
\def\be{\beta}

\def\ka{\kappa}
\def\la{\lambda}

\def\Ga{\Gamma}
\def\Th{\Theta}

\def\CC{{\Bbb C}}
\def\RR{{\Bbb R}}
\def\ZZ{{\Bbb Z}}
\def\FSH{{\cal H}}

{\titlefont \textfont1=\mathtitlefont \textfont2=\mathsymtitlefont
\centerline
{Jacobi polynomials of type $BC$,
Jack polynomials,}
\sPP
\centerline
{limit transitions and $O(\iy)$}}
\bPP
\centerline{Tom H. Koornwinder}
\bLP
\font\smallfont=cmr10 at 10 truept
{\smallfont Extended abstract for the Proceedings of the Conference
``Fourier and Radon transformations on symmetric spaces'' in honor
of Professor S.~Helgason's 65th birthday, Roskilde, Denmark, September 10--12,
1992.}
\bLP
{\bf 1. Limits of Jacobi polynomials.}\quad
Let $\al,\be>-1$. Let $p_n^{(\al,\be)}$ ($n=0,1,2,\ldots\;$) be the
monic orthogonal polynomials \wrt the measure
$(1-x)^\al\,(1+x)^\be\,dx$ on the interval $(-1,1)$ (monic
{\sl Jacobi polynomials}). From the explicit expression as hypergeometric
series we obtain
$$
\lim_{\textstyle{(\al,\be)\to\iy\atop\al/\be\to c}}p_n^{(\al,\be)}(x)=
\left(x+{c-1\over c+1}\right)^n.
$$
In particular,
$$
\lim_{\al\to\iy}p_n^{(\al,\al)}(x)=x^n.
\eqno(1)
$$
Also,
$$
\lim_{\al\to\iy}\al^{\half n}\,p_n^{(\al,\al)}(\al^{-\half n}x)=
h_n(x),
\eqno(2)
$$
where the $h_n$ are the monic orthogonal polynomials \wrt the
measure $e^{-x^2}\,dx$ on $\RR$
(monic {\sl Hermite polynomials}).

\mLP
{\bf 2. A \rep\ of $O(\iy)$.}\quad
Since $p_n^{(\half d -{3\over 2},\half d -{3\over 2})}$ has an interpretation
as spherical function on $O(d)/O(d-1)$, the above limit formulas suggest
that monomials and Hermite polynomials have some interpretation on
$O(\iy)$. This is indeed the case, see for instance McKean [5] and
Matsushima e.a.\ [4]. Let us explain this briefly.
Let $S^{d-1}$ be the sphere of radius $(d/2)^\half$ and midpoint 0 in
$\RR^d$ and consider it as homogeneous space $O(d)/O(d-1)$.
Let $\pi_l^d$ be the irreducible \rep\ of $O(d)$ on the space
$\FSH_l^d$ of spherical harmonics of degree $l$ on $S^{d-1}$.
For $k<d$ put
$$
\FSH_l^{d,k}:=\{f\in\CC[x_1,\ldots,x_k]\mid f|_{S^{d-1}}\in\FSH_l^d\},
$$
i.e.\ the space of $O(d-k)$-invariants in $\FSH_l^d$ (uniquely) extended
to polynomials only depending on the first $k$ coordinates.
This space inherits
an $O(k)$-invariant inner product from $\FSH_l^d$.
As $d\to\iy$, the space $\FSH_l^{d,k}$ tends to a certain $O(k)$-invariant
inner product space $\FSH_l^{\iy,k}$ which is a subspace of
$\CC[x_1,\ldots,x_k]$. Let $\FSH_l^\iy$ be the Hilbert space completion
of $\cup_{k=1}^\iy\FSH_l^{\iy,k}$. Let $O(\iy):=\cup_{d=1}^\iy O(d)$,
considered with the inductive limit topology, where, for $d<d'$, $T\in O(d)$
is embedded in $O(d')$ as $\pmatrix{T&0\cr0&I_{d'-d}}$.
Then $O(\iy)$ acts on $\FSH_l^\iy$
by an irreducible unitary \rep\ $\pi_l^\iy$.

An orthogonal basis of $\FSH_l^\iy$ is given by the polynomials
$$
\{x\mapsto h_{l_1}(x_1)\,h_{l_2}(x_2)\ldots\mid l_1+l_2+\cdots=l\}.
$$
In particular, if $O_r(\iy)$ denotes the subgroup of all $T\in O(\iy)$
which leave the first $r$ standard basis vectors invariant, then
$x\mapsto h_l(x_1)$ is an $O_1(\iy)$-invariant element of $\FSH_l^\iy$,
unique up to a constant factor. This provides a group theoretic interpretation
of (2). The corresponding spherical function can be obtained from the
well-known {\sl addition formula} for Hermite polynomials [3, 10.13 (40)]:
$$
h_l(x_1y_1+\cdots+x_ky_k)=
\sum_{l_1+\cdots l_k=l}{l!\over l_1!\ldots l_k!}\,
y_1^{l_1}\ldots y_k^{l_k}\,h_{l_1}(x_1)\ldots h_{l_k}(x_k),\quad
y_1^2+\cdots y_k^2=1.
$$
Observe that in the above formula the coefficient of $h_l(x_1)$ on the
\RHS\ equals $y_1^l$. This yields the spherical function and, moreover,
a group theoretic interpretation of~(1).

\mLP
{\bf 3. A limit for $BC_n$ type Jacobi polynomials.}\quad
Let $R$ be a root system with Weyl group $W$ and choice of positive roots
$R^+$, and let
$\ka\colon R\to\RR$ be a $W$-invariant mapping. Write
$\rho^{(\ka)}:=\thalf\sum_{\al\in R^+}\ka(\al)\al$.
Let $P^+$ be the set of dominant weights, on which a partial order $\le$ is
induced by $R^+$.
Put $m_\la:=\sum_{\mu\in W.\la}e^\mu$ ($\la\in P^+$).
The {\sl Jacobi polynomial} associated with $R$ (see Beerends \& Opdam [2]
and references given there)
is an element
$p_\la$ ($\la\in P^+$) of the form
$$
p_\la=\sum_{\textstyle{\mu\in P^+\atop\mu\le\la}}\Ga_\mu(\la)\,m_\mu,
\qquad\Ga_\la(\la)=1,
$$
such that it satisfies the differential equation
$$
\bigl({\rm Laplacian}+\sum_{\al\in R^+}\ka(\al)\,\coth(\al/2)\,\pa_\al\bigr)\,
p_\la=\lan\la,\la+2\rho^{(k)}\ran\,p_\la.
$$

For root system $BC_n$ take
$R^+=R_B^+:=\{e_i\}\cup\{2e_i\}\cup\{e_i\pm e_j\}_{i<j}$. Then
$P^+=P_B^+=\{\la\in\ZZ^n\mid\la_1\ge\la_2\ge\ldots\ge\la_n\ge0\}$.
Let $\ka$ send $e_1,2e_1,e_1+e_2$ to $k_1,k_2,k_3$, respectively.
Write $p_\la=p_{B,\la}^{(k_1,k_2,k_3)}$.

For root system $A_{n-1}$ take $R^+=R_A^+=\{e_i-e_j\}_{1\le i<j\le n}$.
Let $\pi$ be the orthogonal projection which maps $\RR^n$ onto the hyperplane
$\{t_1+\cdots t_n=0\}$. Then $P_A^+=\pi(P_B^+)$. Now $\ka$ sends $R$ to
a fixed $k\in\RR$. Write $p_\la=p_{A,\la}^{(k)}$.

We introduce the {\sl Jack polynomial} here by its expression in terms of
$A_{n-1}$-type Jacobi polynomials:
$$
j_\la^{(k)}(t):=
e^{(t_1+\cdots+t_n)(\la_1+\cdots+\la_n)/n}\,
p_{A,\pi(\la)}^{(k)}(\pi(t)),\quad
\la\in P_B^+,\;t\in\RR^n.
$$
The notation used on the \LHS\ for the Jack polynomial,
is not the standard one.
The following limit result generalizes (1).

\mLP
{\bf Theorem 1} (Beerends \& Koornwinder [1]).\quad
Let $-\iy\le a\le1$. Let $\tau,t\in\RR^n$ be related by
$e^{\tau_i}=(a-1)(a-2)^{-1}+(\sinh(\tau_i/2))^2$. Then
$$
\lim_{\textstyle{(k_1,k_2)\to\iy\atop k_1/k_2\to-a}}
p_{B,\la}^{(k_1,k_2,k_3)}(t)=\const j_\la^{(k_3)}(\tau).
$$

\mLP
{\bf 4. An interpretation of Jack polynomials on $O(\iy)$.}\quad
We first summarize some results of Olshanski{\u\i} [6].
Denote the left upper $r\times r$ submatrix of the matrix $T\in O(\iy)$
by $\Th_r(T)$. Then the mapping which sends a double coset
$O_r(\iy).T.O_r(\iy)$ to $\Th_r(T)$ is well-defined and gives a bijection
of the double coset space $O(\iy)//O_r(\iy)$ onto the semigroup $\Ga(r)$
of all $r\times r$ matrices with norm $\le1$.
Let $\pi$ be an irreducible unitary \rep\ of $O(\iy)$ on a Hilbert space
$\FSH$ \st\ the subspace $\FSH_r$ of $O_r(\iy)$-fixed elements in $\FSH$
contains nonzero elements.
Let $P_r$ be the orthogonal projection of $\FSH$ onto $\FSH_r$.

\mLP
{\bf Theorem 2} (Olshanski{\u\i}).\quad
The mapping
$\Th_r(T)\mapsto P_r\pi(T)|_{\FSH_r}$ ($T\in O(\iy)$) is a semigroup
\rep\ of $\Ga(r)$ on $\FSH_r$.
Moreover, the space $\FSH_r$ has finite dimension and the semigroup \rep\
has a unique extension to an
irreducible group \rep\ of $GL(n,\RR)$ on $\FSH_r$.

\mPP
Suppose now that the space of $O(r)\times O_r(\iy)$-invariant elements
of $\FSH$ has dimension one. Then we can speak about the
corresponding spherical function on $O(\iy)$ \wrt $O(r)\times O_r(\iy)$.
By the above Theorem this spherical function can be immediately
connected with a spherical function on $GL(r,\RR)$ \wrt $O(r)$
corresponding to a finite-dimensional irreducible \rep\ of $GL(r,\RR)$.
These last spherical functions were called {\sl zonal polynomials} by
A.~T. James and they are special instances of Jack polynomials.

For certain values of $k_1,k_2,k_3$ the $BC_r$-type Jacobi polynomial
can be interpreted as a spherical function on $O(d)$ \wrt
$O(r)\times O(d-r)$. It is plausible to suppose that such a spherical
function tends, for $d\to\iy$, to a certain spherical function on
$O(\iy)$ \wrt $O(r)\times O_r(\iy)$. When we connect this, in its turn,
to a spherical function on $GL(r,\RR)$ \wrt $O(r)$ and write this as a
Jack polynomial then we would arrive at a group theoretic interpretation
of Theorem 1 in special cases.

\mLP
{\bf References}\frenchspacing

\sLP
[1] R. J. Beerends \& T. H. Koornwinder,
{\sl Analysis on root systems: $A_{n-1}$ as limit case of $BC_n$},
in preparation.

\sLP
[2] R. J. Beerends \& E. M. Opdam,
{\sl Certain hypergeometric series related to the root system $BC$},
Trans. Amer. Math. Soc., to appear.

\sLP
[3] A. Erd{\'e}lyi,
W. Magnus,
F. Oberhettinger \&
F. G. Tricomi,
{\sl Higher transcendental functions, Vol. II},
McGraw-Hill, 1953.

\sLP
[4] H. Matsushima,
K. Okamoto \&
T. Sakurai,
{\sl On a certain class of irreducible unitary \rep s of the infinite
dimensional rotation group I},
Hiroshima Math. J. 11 (1981), 181--193.

\sLP
[5] H. P. McKean,
{\sl Geometry of differential space},
Ann. Probab. 1 (1973), 197--206.

\sLP
[6] G.I. Olshanski{\u\i},
{\sl Representations of infinite-dimensional classical groups,
limits of enveloping algebras and Yangians},
in: {\sl Topics in representation theory},
A. A. Kirillov (ed.),
Advances in Soviet Mathematics 2, 1991.

\bPP
\obeylines\parindent 5truecm
{\addressfont
University of Amsterdam, Faculty of Mathematics and Computer Science
Plantage Muidergracht 24, 1018 TV Amsterdam, The Netherlands
email {\ttaddressfont thk@fwi.uva.nl}}

\bye